\documentclass{amsart}%
\usepackage{graphicx}
\usepackage{amscd}
\usepackage{amsmath}
\usepackage{amsfonts}
\usepackage{amssymb}%
\newtheorem{theorem}{Theorem}
\theoremstyle{plain}

\newtheorem{corollary}{Corollary}

\newtheorem{definition}{Definition}
\newtheorem{example}{Example}

\newtheorem{lemma}{Lemma}

\newtheorem{proposition}{Proposition}

\numberwithin{equation}{section}

\begin{document}
\title[Almost summing mappings]{Almost summing mappings}
\author{Daniel Pellegrino}
\address[Daniel Pellegrino]{Depto de Matem\'{a}tica e Estat\'{\i}stica- Caixa Postal
10044- UFPB Campus II- Campina Grande-PB-Brasil}
\email{dmp@dme.ufcg.edu.br} 

\begin{abstract}
We introduce a general definition of almost $p$-summing mappings and give
several concrete examples of such mappings. Some known results are
considerably generalized and we present various situations in which the space
of almost $p$-summing multilinear mappings coincides with the whole space of
continuous multilinear mappings.

\end{abstract}
\maketitle

\section{Introduction}

The rapid development of the theory of absolutely summing linear mappings has
lead to the study of innumerous new classes of multilinear mappings and
polynomials between Banach spaces (see \cite{Matos2},\cite{MatosJ}%
,\cite{BOtelho25},\cite{Bombal}). Recently, Botelho \cite{BOtelho25} and
Botelho-Braunss-Junek \cite{Junek} introduced the concept of almost
$p$-summing multilinear mappings and gave the first examples and properties of
such mappings. The recent work of Matos \cite{Matos}, concerning absolutely
summing arbitrary mappings, turns natural to ask whether it is possible to
follow the same line of thought with almost $p$-summing mappings. In this
paper we will present a more general definition of almost $p$-summing
mappings, several examples and a natural version of a Dvoretzky-Rogers Theorem
for this kind of applications. It will be shown that almost $p$-summing
multilinear mappings are much more common than it was known until now. For
example, we prove that every continuous $n$-linear mapping from $C(K)\times
...\times C(K)$ into a Banach space $F$ is almost $2$-summing, generalizing a
recent result obtained in \cite{Junek}. This paper also analyzes the
connections of almost $p$-summing mappings and type/cotype and provides
various examples of analytic almost $p$-summing mappings.

\section{Absolutely summing mappings}

Throughout this paper $E,E_{1},...,E_{n},F$ will stand for Banach spaces. For
$p\in\lbrack1,\infty\lbrack,$ the linear space of all sequences $(x_{j}%
)_{j=1}^{\infty}$ in $E$ such that
\[
\Vert(x_{j})_{j=1}^{\infty}\Vert_{p}=(\sum_{j=1}^{\infty}\Vert x_{j}\Vert
^{p})^{\frac{1}{p}}<\infty
\]
will be denoted by $l_{p}(E).$ We will denote by $l_{p}^{w}(E)$ the linear
space formed by the sequences $(x_{j})_{j=1}^{\infty}$ in $E$ such that
$(<\varphi,x_{j}>)_{j=1}^{\infty}\in l_{p}(\mathbb{K}),$ for every continuous
linear functional $\varphi:E\rightarrow\mathbb{K}.$ We also define
$\Vert.\Vert_{w,p}$ in $l_{p}^{w}(E)$ by
\[
\Vert(x_{j})_{j=1}^{\infty}\Vert_{w,p}:=\sup_{\varphi\in B_{E^{\prime}}}%
(\sum_{j=1}^{\infty}\mid<\varphi,x_{j}>\mid^{p})^{\frac{1}{p}}.
\]
The linear subspace of $l_{p}^{w}(E)$ of all sequences $(x_{j})_{j=1}^{\infty
}\in l_{p}^{w}(E),$ such that
\[
\lim_{m\rightarrow\infty}\Vert(x_{j})_{j=m}^{\infty}\Vert_{w,p}=0
\]
will be denoted by $l_{p}^{u}(E)$. The sequences in $l_{p}^{u}(E)$ are called
unconditionally $p$-summable.

The multilinear theory of absolutely summing mappings was first sketched by
Pietsch in \cite{Pietsch} and has been broadly explored (see \cite{Tonge},
\cite{Matos2}, \cite{Floret}). The next definition can be found in
\cite{Matos2}.

\begin{definition}
\label{deff} A multilinear mapping $T:E_{1}\times...\times E_{n}\rightarrow F$
is absolutely\emph{\ }$(p;q_{1},...,q_{n})$-summing if
\[
(T(x_{j}^{(1)},...,x_{j}^{(n)}))_{j=1}^{\infty}\in l_{p}(F)
\]
for every $(x_{j}^{(s)})_{j=1}^{\infty}\in l_{q_{s}}^{w}(E_{s})$, $s=1,...,n.$
An $n$-homogeneous polynomial $P:E\rightarrow F$ is absolutely $(p;q)$-summing
if
\[
(P(x_{j}))_{j=1}^{\infty}\in l_{p}(F)
\]
whenever $(x_{j})_{j=1}^{\infty}\in l_{q}^{w}(E).$
\end{definition}

It is worth observing that, in Definition \ref{deff}, there is no difference
if we consider $l_{q_{s}}^{u}(E)$ $(l_{q}^{u}(E))$ instead of $l_{q_{s}}%
^{w}(E)$ $(l_{q}^{w}(E))$ (see \cite[Proposition 2.4]{Matos2} for polynomials,
and the multilinear case is analogous).

The following well known characterization can be found in \cite[Theorem
1.2(ii)]{Botelho2}, and is sometimes useful.

\begin{theorem}
Let $T:E_{1}\times...\times E_{n}\rightarrow F$ be a multilinear mapping. The
following statements are equivalent:

(1) $T$ is absolutely $(p;q_{1},...,q_{n})$-summing.

(2) There exists $L>0$ such that for every natural $k$ and any $x_{j}^{(l)}\in
E_{l},$%
\begin{equation}
(\sum_{j=1}^{k}\Vert T(x_{j}^{(1)},...,x_{j}^{(n)})\Vert^{p})^{\frac{1}{p}%
}\leq L\Vert(x_{j}^{(1)})_{j=1}^{k}\Vert_{w,q_{1}}...\Vert(x_{j}^{(n)}%
)_{j=1}^{k}\Vert_{w,q_{n}}. \label{aqqq}%
\end{equation}

The least $L>0$ for which inequality (\ref{aqqq}) always holds defines a norm
for the space of absolutely $(p;q_{1},...,q_{n})$-summing multilinear
mappings. This norm will be denoted by $\Vert.\Vert_{as(p;q)}.$ A
characterization for $n$-homogeneous polynomials is analogous.
\end{theorem}

Inspired on the work of Matos \cite{Matos1}, we introduce the following
concept, which generalizes Definition \ref{deff}, as we will see later.

\begin{definition}
\label{n}An arbitrary mapping $f$:$E\rightarrow F$ is absolutely
$(p,q)$-summing at $a$ if there exist $M_{a}>0,$ $\delta_{a}>0$ and $r_{a}>0$
so that
\[
\sum_{j=1}^{k}\Vert f(a+x_{j})-f(a)\Vert^{p}\leq M_{a}\Vert(x_{j})_{j=1}%
^{k}\Vert_{w,q}^{r_{a}}%
\]
for all $k$ and $\Vert(x_{j})_{j=1}^{k}\Vert_{w,q}<\delta_{a}.$
\end{definition}

\begin{theorem}
\label{sh}If $F$ has cotype $q$, $E$ is an $\mathcal{L}_{\infty,\lambda}$
space and $f:E\rightarrow F$ is analytic at $a$, then $f$ is absolutely
$(q;2)$-summing at $a$.
\end{theorem}

Proof. Since $f$ is analytic at $a$, there are $C\geq0$ and $c>0$ such that\
\[
\Vert\frac{1}{k!}\overset{\wedge}{d^{k}}f(a)\Vert\leq Cc^{k}\text{ for every
}k.
\]
A recent result of D. Perez (see \cite{Perez}) states that whenever each
$E_{j}$ is an $\mathcal{L}_{\infty,\lambda_{j}}$ space, every continuous
$n$-linear ($n\geq2$) mapping $T$, from $E_{1}\times...\times E_{n}$ into
$\mathbb{K,}$ is absolutely $(1;2,...,2)$-summing and
\begin{equation}
\Vert T\Vert_{as(1;2,...,2)}\leq K_{G}3^{\frac{n-2}{2}}\Vert T\Vert
\prod\limits_{j=1}^{n}\lambda_{j}. \label{DPerez}%
\end{equation}
Using the polynomial version of this result, it is not hard to prove that (see
\cite[Theorem 4]{Pellegrino}) whenever $F$ has finite cotype $q$, every
bounded $n$-homogeneous ($n\geq2$) polynomial $P:E\rightarrow F$ is absolutely
$(q;2)$-summing and $\Vert P\Vert_{as(q;2)}\leq C_{q}(F)K_{G}3^{\frac{n-2}{2}%
}\Vert P\Vert\lambda^{n},$ where $C_{q}(F)$ and $K_{G}$ are the cotype's
constant of $F$ and Grothendieck's constant, respectively.

For $n=1$, we still have $\mathcal{L}(E;F)=\mathcal{L}_{as(q;2)}(E;F),$ which
is a particular case of a result due to Dubinsky-Pe\l czy\'{n}ski-Rosenthal
(case $q=2$) and Maurey (case $q>2$) (see \cite[Theorem 11.14 (a) and (b)
]{Diestel}). So, for every natural $n$, there exist positive $D$ and $D_{1}$
so that
\[
\Vert\frac{1}{k!}\overset{\wedge}{d^{k}}f(a)\Vert_{as(q;2)}\leq D_{1}%
D^{k}\Vert\frac{1}{k!}\overset{\wedge}{d^{k}}f(a)\Vert.
\]
Hence, if $\delta_{a}$ is the radius of convergence of $f$ around $a$, then,
whenever $(x_{j})_{j=1}^{m}$ is such that $\Vert(x_{j})_{j=1}^{m}\Vert
_{w,1}\leq\min\{\frac{1}{2D},\delta_{a}\},$ we have
\begin{align*}
(\sum_{j=1}^{m}\Vert f(a+x_{j})-f(a)\Vert^{q})^{\frac{1}{q}}  &  =\sum
_{j=1}^{m}(\Vert\sum_{k=1}^{\infty}\frac{1}{k!}\overset{\wedge}{d^{k}%
}f(a)(x_{j})\Vert^{q})^{\frac{1}{q}}\\
&  \leq\sum_{k=1}^{\infty}[\sum_{j=1}^{m}\Vert\frac{1}{k!}\overset{\wedge
}{d^{k}}f(a)(x_{j})\Vert^{q}]^{\frac{1}{q}}\\
&  \leq\sum_{k=1}^{\infty}\Vert\frac{1}{k!}\overset{\wedge}{d^{k}}%
f(a)\Vert_{as(q;2)}\Vert(x_{j})_{j=1}^{m}\Vert_{w,2}^{k}\\
&  \leq D_{1}\Vert(x_{j})_{j=1}^{m}\Vert_{w,2}\sum_{k=1}^{\infty}\frac{D^{k}%
}{2^{k-1}D^{k-1}}=2DD_{1}\Vert(x_{j})_{j=1}^{m}\Vert_{w,2}.
\end{align*}
$\Box$

Several other results concerning absolutely summing analytic mappings can be
found in \cite{Floret} and \cite{Pellegrino}.

\begin{proposition}
\label{y}If $f:E\rightarrow F$ is absolutely $(p;q)$-summing at $a$, then $f$
is so that $(f(a+x_{j})-f(a))_{j=1}^{\infty}$ $\in l_{p}(F)$ whenever
$(x_{j})_{j=1}^{\infty}$ is unconditionally $q$-summable.
\end{proposition}

Proof. Let $f$ be $(p;q)$-summing at $a$. For any $(x_{j})_{j=1}^{\infty}\in
l_{p}^{u}(E)$, we have
\[
\lim_{k,m\rightarrow\infty}(\sum_{j=k}^{m}\Vert f((a+x_{j})-f(a)\Vert
^{p})^{\frac{1}{p}}\leq\lim_{k,m\rightarrow\infty}C_{a}\Vert(x_{j})_{j=k}%
^{m}\Vert_{w,p}^{r_{a}}=0
\]
and, by the completeness of $l_{p}(F),$ we obtain $(f(a+x_{j})-f(a))_{j=1}%
^{\infty}\in l_{p}(F)$.$\Box$

An immediate outcome of Proposition \ref{y} is that Definition \ref{n} applied
for $n$-homogeneous polynomials and the usual definition of absolutely
$(p,q)$-summing polynomials coincides at $a=0$. In order to prove that
Definition \ref{n} for $n$-linear mappings actually generalizes the standard
definition (Definition \ref{deff}) of absolutely $(p;q_{1},...,q_{n})$-summing
multilinear mappings for $q_{1}=...=q_{n}=q$, we need the following Lemma,
which is a simple consequence of the Open Mapping Theorem.

\begin{lemma}
\label{1}$l_{q}^{u}(E_{1}\times...\times E_{n})$ is isomorphic to $l_{q}%
^{u}(E_{1})\times....\times l_{q}^{u}(E_{n}).$
\end{lemma}

\begin{proposition}
\label{nn}An $n$-linear mapping $T$ is $(p;q,...,q)$-summing in the usual
sense if, and only if, it is absolutely $(p;q)$-summing at the origin in the
sense of Definition \ref{n}.
\end{proposition}

Proof. Consider an absolutely $(p;q)$-summing (in the sense of Definition
\ref{n}, at the origin) $n$-linear mapping, $T:E_{1}\times...\times
E_{n}\rightarrow F$. Then, given $(x_{j}^{(1)})_{j=1}^{\infty}\in l_{q}%
^{u}(E_{1}),....,(x_{j}^{(n)})_{j=1}^{\infty}\in l_{q}^{u}(E_{n}),$ we have
$(x_{j}^{(1)},...,x_{j}^{(n)})_{j=1}^{\infty}\in l_{q}^{u}(E_{1}%
\times....\times E_{n}).$ Hence, by Proposition \ref{y}, $(T(x_{j}%
^{(1)},....,x_{j}^{(n)}))_{j=1}^{\infty}\in l_{p}(F).$ Thus, by the usual
definition, it follows that $T$ is absolutely $(p;q,...,q)$-summing .

Conversely, consider an absolutely $(p;q,...,q)$-summing $n$-linear mapping
$T$ in the usual meaning. Then, if $x_{1}^{(l)},...,x_{k}^{(l)}\in
E_{l},l=1,...,n,$ we have
\[
(\sum_{j=1}^{k}\Vert(T(x_{j}^{(1)},...,x_{j}^{(n)})\Vert^{p})^{\frac{1}{p}%
}\leq C\Vert(x_{j}^{(1)})_{j=1}^{k}\Vert_{w,q}....\Vert(x_{j}^{(n)})_{j=1}%
^{k}\Vert_{w,q}.
\]
Therefore, since $l_{q}^{u}(E_{1}\times....\times E_{n})$ is isomorphic to
$l_{q}^{u}(E_{1})\times....\times l_{q}^{u}(E_{n}),$ it follows that there
exists $C_{1}>0$ so that, for every $k$,
\[
\Vert(x_{j}^{(1)},...,x_{j}^{(n)})_{j=1}^{k}\Vert_{w,q}\geq C_{1}(\Vert
(x_{j}^{(1)})_{j=1}^{k}\Vert_{w,q}+...+\Vert(x_{j}^{(n)})_{j=1}^{k}\Vert
_{w,q})
\]
and
\begin{align*}
\Vert(x_{j}^{(1)},...,x_{j}^{(n)})_{j=1}^{k}\Vert_{w,q}^{n}  &  \geq C_{1}%
^{n}(\Vert(x_{j}^{(1)})_{j=1}^{k}\Vert_{w,q}+...+\Vert(x_{j}^{(n)})_{j=1}%
^{k}\Vert_{w,q})^{n}\\
&  \geq C_{1}^{n}(\Vert(x_{j}^{(1)})_{j=1}^{k}\Vert_{w,q}...\Vert(x_{j}%
^{(n)})_{j=1}^{k}\Vert_{w,q})\\
&  \geq\frac{C_{1}^{n}}{C}(\sum_{j=1}^{k}\Vert(T(x_{j}^{(1)},...,x_{j}%
^{(n)})\Vert^{p})^{\frac{1}{p}}%
\end{align*}
and so $T$ is absolutely $(p;q)$-summing in the sense of Definition \ref{n}.
$\Box$

\section{Almost summing mappings}

Considering the Rademacher functions $(r_{j}(t))_{j=1}^{\infty},$ we say that
the sequence $(x_{j})_{j=1}^{\infty}$ of points of $E$ is almost
unconditionally summable if $\sum\limits_{j=1}^{\infty}r_{j}(t)x_{j}\in
L_{p}([0,1],E)$ for some$,$ and then for all $p$, $0<p<\infty.$

\begin{definition}
\label{esperar}(Botelho \cite{BOtelho25}) An $n$-linear mapping $T:E_{1}%
\times...\times E_{n}\rightarrow F$ is said to be almost $(p_{1},...,p_{n}%
)$-summing if there exists $C\geq0$ such that
\[
(\int\limits_{0}^{1}\Vert\sum_{j=1}^{k}T(x_{j}^{(1)},...,x_{j}^{(n)}%
)r_{j}(t)\Vert^{2}dt)^{\frac{1}{2}}\leq C\Vert(x_{j}^{(1)})_{j=1}^{k}%
\Vert_{w,p_{1}}...\Vert(x_{j}^{(n)})_{j=1}^{k}\Vert_{w,p_{n}}%
\]
for every $k$ and any $x_{j}^{(l)}$ in $E_{l},l=1,...,n$ and $j=1,...,k.$ An
$n$-homogeneous polynomial $P:E\rightarrow F$ is said almost $p$-summing when
$\overset{\vee}{P}$ is almost $(p,...,p)$-summing. The space of all almost
$p$-summing polynomials is denoted by $\mathcal{P}_{al,p}(^{n}E;F).$
\end{definition}

\begin{theorem}
\label{arch}(\cite[Theorem 3.3]{Junek})For $1\leq p\leq2n$ and $P\in
\mathcal{P}_{al,p}(^{n}E;F),$ the following properties are equivalent:

(i) $P$ is almost $p$-summing.

(ii) $P$ maps unconditionally $p$-summable sequences in $E$ into almost
unconditionally summable sequences in $F$.
\end{theorem}

The following definition is a natural generalization of Definition
\ref{esperar} and allows us to give examples of analytic almost $p$-summing mappings.

\begin{definition}
\label{RPM2}A mapping $f:E\rightarrow F$ is said to be almost $p$-summing at
$a\in E$ if there exist $C_{a}>0$, $\epsilon_{a}>0$ and $r_{a}>0$ \ such that
\[
(\int\limits_{0}^{1}\Vert\sum_{j=1}^{k}(f(a+x_{j})-f(a))r_{j}(t)\Vert
^{2}dt)^{\frac{1}{2}}\leq C_{a}\Vert(x_{j})_{j=1}^{k}\Vert_{w,p}^{r_{a}}%
\]
for every natural $k$, any $x_{1},...,x_{k}$ in $E$ and $\Vert(x_{j}%
)_{j=1}^{k}\Vert_{w,p}<\epsilon_{a}.$ If $f$ is almost $p$-summing at every
$a\in E$, \ we say that $f$ is almost $p$-summing everywhere$.$
\end{definition}

It is worth observing that if $f$ is almost $p$-summing at $a$, then $f$ is
continuous at $a.$ The space of all polynomials from $E$ into $F$ which are
almost $p$-summing everywhere will be denoted by $\mathcal{P}_{al,p(E)}%
(^{n}E;F).$

\begin{proposition}
\label{Flamante}If $f:E\rightarrow F$ is \ almost $p$-summing at $a$, then $f$
is so that $(f(a+x_{j})-f(a))_{j=1}^{\infty}$ is almost unconditionally
summable whenever $(x_{j})_{j=1}^{\infty}$ is unconditionally $p$-summable.
\end{proposition}

Proof. Analogous to the proof of Proposition \ref{y}.

An immediate outcome of Theorem \ref{arch} and Proposition \ref{Flamante} is
that Definitions \ref{RPM2} and \ref{esperar} coincides for $n$-homogeneous
polynomials and $a=0$. The proof that Definition \ref{RPM2}, for $a=0,$
generalizes Definition \ref{esperar}, for multilinear mappings and
$p_{1}=...=p_{n}=p$, is similar to the proof of Proposition \ref{nn}.

\begin{proposition}
If $P\in\mathcal{P}(^{n}E;F)$, then $P\in\mathcal{P}_{al,p(E)}(^{n}%
E;F)\Leftrightarrow\overset{\vee}{P}\in\mathcal{L}_{al,p(E)}(^{n}E,F).$
\end{proposition}

Proof. Suppose that $P\in\mathcal{P}_{al,p(E)}(^{n}E;F).$ Then, by the
polarization formula,
\[
\overset{\vee}{P}(a_{1}+x_{j}^{(1)},...,a_{n}+x_{j}^{(n)})-\overset{\vee}%
{P}(a_{1},...,a_{n})=
\]%
\begin{align*}
&  =[\frac{1}{n!2^{n}}\sum_{e_{i}=1,-1}e_{1}...e_{n}P(e_{1}(a_{1}+x_{j}%
^{(1)})+...+e_{n}(a_{n}+x_{j}^{(n)})]-\\
&  -[\frac{1}{n!2^{n}}\sum_{e_{i}=1,-1}e_{1}...e_{n}P(e_{1}a_{1}%
+...+e_{n}a_{n})]\\
&  =\frac{1}{n!2^{n}}\sum_{e_{i}=1,-1}e_{1}...e_{n}[P((e_{1}a_{1}%
+...+e_{n}a_{n})+(e_{1}x_{j}^{(1)}+...+e_{n}x_{j}^{(n)}))-\\
&  -P(e_{1}a_{1}+...+e_{n}a_{n})].
\end{align*}
For any $(x_{j}^{(1)})_{j=1}^{k},...,(x_{j}^{(n)})_{j=1}^{k}$ , in order to
simplify notation, we will write
\[
A=(\int\limits_{0}^{1}\Vert\sum_{j=1}^{k}(\overset{\vee}{P}(a_{1}+x_{j}%
^{(1)},...,a_{n}+x_{j}^{(n)})-\overset{\vee}{P}(a_{1},...,a_{n}))r_{j}%
(t)\Vert^{2}dt)^{\frac{1}{2}}.
\]
Lemma \ref{1} asserts that there exists $D>0$ so that
\[
(\Vert(x_{j}^{(1)})_{j=1}^{k}\Vert_{w,p}+...+\Vert(x_{j}^{(n)})_{j=1}^{k}%
\Vert_{w,p})\leq D\Vert(x_{j}^{(1)},...,x_{j}^{(n)})_{j=1}^{k}\Vert_{w,p}%
\]
for every $k.$ Now suppose
\[
\Vert(x_{j}^{(1)},...,x_{j}^{(n)})_{j=1}^{k}\Vert_{w,p}<\frac{1}{D}\min
_{e_{i}=-1,1}\{1,\epsilon_{e_{1}a_{1}+...+e_{n}a_{n}}\},
\]
where the $\epsilon_{e_{1}a_{1}+...+e_{n}a_{n}}$ are given by Definition
\ref{RPM2} applied to $P$. Then, for any choice of $-1$ and $1$ for $e_{j},$
we have
\[
\Vert(e_{1}x_{j}^{(1)}+...+e_{n}x_{j}^{(n)})_{j=1}^{k}\Vert_{w,p}<\min
_{e_{i}=-1,1}\{1,\epsilon_{e_{1}a_{1}+...+e_{n}a_{n}}\}.
\]
Therefore,
\begin{align*}
A  &  =(\int\limits_{0}^{1}\Vert\sum_{j=1}^{k}\frac{1}{n!2^{n}}\sum
_{e_{i}=1,-1}e_{1}...e_{n}[P((e_{1}a_{1}+...e_{n}a_{n})+(e_{1}x_{j}%
^{(1)}+...+e_{n}x_{j}^{(n)}))-\\
&  -P(e_{1}a_{1}+...+e_{n}a_{n})]r_{j}(t)\Vert^{2}dt)^{\frac{1}{2}}\\
&  \leq\frac{1}{n!2^{n}}\sum_{e_{i}=1,-1}(\int\limits_{0}^{1}\Vert\sum
_{j=1}^{k}e_{1}...e_{n}[P((e_{1}a_{1}+...+e_{n}a_{n})+(e_{1}x_{j}%
^{(1)}+...+e_{n}x_{j}^{(n)}))-\\
&  -P(e_{1}a_{1}+...+e_{n}a_{n})]r_{j}(t)\Vert^{2}dt)^{\frac{1}{2}}\\
&  \leq\frac{1}{n!2^{n}}\sum_{e_{i}=1,-1}C_{e_{1}a_{1}+...+e_{n}a_{n}}%
\Vert(e_{1}x_{j}^{(1)}+...+e_{n}x_{j}^{(n)})_{j=1}^{k}\Vert_{w,p}%
^{r_{(e_{1}a_{1}+...+e_{n}a_{n})}}\\
&  \leq\frac{1}{n!2^{n}}\sum_{e_{i}=1,-1}C_{e_{1}a_{1}+...+e_{n}a_{n}}%
(\Vert(x_{j}^{(1)})_{j=1}^{k}\Vert_{w,p}+...+\Vert(x_{j}^{(n)})_{j=1}^{k}%
\Vert_{w,p})^{r_{(e_{1}a_{1}+...+e_{n}a_{n})}}\\
&  \leq\frac{1}{n!2^{n}}\sum_{e_{i}=1,-1}C_{e_{1}a_{1}+...+e_{n}a_{n}%
}D^{_{r_{_{(e_{1}a_{1}+...+e_{n}a_{n})}}}}\Vert(x_{j}^{(1)},...,x_{j}%
^{(n)})_{j=1}^{k}\Vert_{w,p}^{^{r_{(e_{1}a_{1}+...+e_{n}a_{n})}}}\\
&  \leq D_{1}\Vert(x_{j}^{(1)},...,x_{j}^{(n)})_{j=1}^{k}\Vert_{w,p}%
^{\min\{r_{(e_{1}a_{1}+...+e_{n}a_{n})}\}}%
\end{align*}
if $\Vert(x_{j}^{(1)},...,x_{j}^{(n)})_{j=1}^{k}\Vert_{w,p}<\delta=\frac{1}%
{D}\min_{e_{i}=-1,1}\{1,\epsilon_{e_{1}a_{1}+...+e_{n}a_{n}}\}.$The converse
is obvious. $\Box$

Naturally, the concepts of type and cotype give us the next Proposition.

\begin{proposition}
\label{brass} If $F$ has type $q$, then every absolutely $(q;p)$-summing
mapping (at $a$) is almost $p$-summing at $a.$ On the other hand, if $F$ has
finite cotype $r$, then every almost $p$-summing mapping (at $a$) is
$(r;p)$-summing at $a$.
\end{proposition}

\begin{corollary}
If $F$ is a Hilbert space and $E$ is an $\mathcal{L}_{\infty}$ space, then
every $f:E\rightarrow F,$ analytic at $a$, is almost $2$-summing at $a$. In
particular, under the same hypothesis, every entire mapping $f:E\rightarrow F$
is almost $2$-summing everywhere.
\end{corollary}

Proof. Since $\cot F=2,$ by Theorem \ref{sh}, $\ f$ is absolutely
$(2;2)$-summing at $a.$ Besides, since $F$ has type $2$, then $f$ is almost
$2$-summing at $a$, by Proposition \ref{brass}.$\Box$

In order to give the other examples of analytic almost summing mappings, the
next Proposition will be useful.

\begin{proposition}
\label{Penelope}If $f$ is such that there exist $C,$ $\delta,r>0$ so that
\[
\Vert(f(a+x_{j})-f(a))_{j=1}^{k}\Vert_{w,1}\leq C\Vert(x_{j})_{j=1}^{k}%
\Vert_{w,p}^{r}%
\]
for any natural $k,$ every $x_{1},...,x_{k}$ in $E$ and $\Vert(x_{j}%
)_{j=1}^{k}\Vert_{w,p}<\delta,$ then $f$ is almost $p$-summing at $a$.
\end{proposition}

Proof.
\[
(\int\limits_{0}^{1}\Vert\sum_{j=1}^{k}(f(a+x_{j})-f(a))r_{j}(t)\Vert
^{2}dt)^{\frac{1}{2}}\leq\sup_{t\in\lbrack0,1]}\Vert\sum_{j=1}^{k}%
(f(a+x_{j})-f(a))r_{j}(t)\Vert=
\]%
\begin{align*}
&  =\sup_{t\in\lbrack0,1]}\sup_{\varphi\in B_{E}%
\acute{}%
}\mid<\varphi,\sum_{j=1}^{k}(f(a+x_{j})-f(a))r_{j}(t)>\mid\\
&  \leq\Vert(f(a+x_{j})-f(a))_{j=1}^{k}\Vert_{w,1}\leq C\Vert(x_{j})_{j=1}%
^{k}\Vert_{w,p}^{r}%
\end{align*}
for $\Vert(x_{j})_{j=1}^{k}\Vert_{w,p}<\delta.\Box$

In \cite[Corollary 6.3]{BOtelho25} it is stated that regardless of the
positive integer $n$, every absolutely $(1;2)$-summing $n$-homogeneous
polynomial is almost $2$-summing. It is worth remarking that, when $f$ is a
polynomial, $a=0$ and $p=2$, Proposition \ref{Penelope} is a significant
improvement of \cite[Corollary 6.3]{BOtelho25}, since in Proposition
\ref{Penelope} we just need a weak estimate whereas in \cite[Corollary
6.3]{BOtelho25} we need a norm estimate. As we will see later in Corollary
\ref{finno}, the aforementioned Proposition is the key of innumerous new
Coincidence Theorems which will generalize the few Coincidence Theorems known
until now (see \cite[Proposition 7.1]{BOtelho25},\cite[Proposition 5.1]%
{Junek}). The next Corollary give other examples of almost $p$-summing
analytic mappings.

\begin{corollary}
Let $E$ be an $\mathcal{L}_{\infty,\lambda}$ space and $F$ be an arbitrary
Banach space. Every mapping $g:E\rightarrow F,$ analytic at $a$, such that
$dg(a)=0$ is almost $2$-summing at $a$.
\end{corollary}

Proof. Let $C$ and $c$ be such that
\[
\Vert\frac{1}{k!}\overset{\wedge}{d^{k}}g(a)\Vert\leq Cc^{k}\text{ for every
}k\geq1.
\]
Then, for any bounded linear functional $\varphi,$ defined on $F,$ we obtain
\[
\Vert\frac{1}{k!}\overset{\wedge}{d^{k}}\varphi g(a)\Vert=\Vert\varphi\frac
{1}{k!}\overset{\wedge}{d^{k}}g(a)\Vert\leq Cc^{k}\Vert\varphi\Vert\text{ for
every }k\geq1.
\]
By (\ref{DPerez}) we have
\[
\Vert\frac{1}{k!}\overset{\wedge}{d^{k}}\varphi g(a)\Vert_{as(1;2)}\leq
K_{G}3^{\frac{k-2}{2}}\lambda^{k}Cc^{k}\Vert\varphi\Vert\text{ for every
}k\geq2.
\]
Therefore, defining $\delta_{a}$ as the radius of convergence of $g$ around
$a$, if we assume $(x_{j})_{j=1}^{m}$ such that
\[
\Vert(x_{j})_{j=1}^{m}\Vert_{w,2}\leq\delta=\min\{\frac{1}{(2\sqrt{3}\lambda
c)},\delta_{a}\},
\]
we obtain
\begin{align*}
\sum_{j=1}^{m}  &  \mid\varphi g(a+x_{j})-\varphi g(a)\mid\leq\sum
_{k=2}^{\infty}\Vert\frac{1}{k!}\overset{\wedge}{d^{k}}\varphi g(a)\Vert
_{as(1;2)}\Vert(x_{j})_{j=1}^{m}\Vert_{w,2}^{k}\\
&  =\Vert(x_{j})_{j=1}^{m}\Vert_{w,2}\sum_{k=2}^{\infty}\Vert\frac{1}%
{k!}\overset{\wedge}{d^{k}}\varphi g(a)\Vert_{as(1;2)}\Vert(x_{j})_{j=1}%
^{m}\Vert_{w,2}^{k-1}\\
&  \leq\Vert(x_{j})_{j=1}^{m}\Vert_{w,2}\sum_{k=2}^{\infty}\frac{K_{G}%
3^{\frac{k-2}{2}}\lambda^{k}Cc^{k}\Vert\varphi\Vert}{(2\sqrt{3}\lambda
c)^{k-1}}\leq D\Vert(x_{j})_{j=1}^{m}\Vert_{w,2}%
\end{align*}
for every $\varphi\in B_{F}%
\acute{}%
$ and every $m$. Therefore,
\[
\Vert(g(a+x_{j})-g(a))_{j=1}^{m}\Vert_{w,1}\leq D\Vert(x_{j})_{j=1}^{m}%
\Vert_{w,2}\text{ }%
\]
regardless of the $\Vert(x_{j})_{j=1}^{m}\Vert_{w,p}<\delta$ , and
$x_{1},...,x_{m}.$ Now, Proposition \ref{Penelope} yields the result. $\Box$

In \cite[Proposition 5.1]{Junek} it is shown that if $E$ is an $\mathcal{L}%
_{\infty}$ space then $\mathcal{L}(^{2}E;\mathbb{K})=\mathcal{L}_{al,2}%
(^{2}E;\mathbb{K}).$ Next corollary shows that the aforementioned result is
still valid for vector valued $n$-linear mappings, for every $n\geq2$.

\begin{corollary}
\label{finno} If $E$ is an $\mathcal{L}_{\infty}$ space and $n\geq2$, then for
every Banach space $F$ we have
\begin{equation}
\mathcal{P}_{al,2}(^{n}E;F)=\mathcal{P}(^{n}E;F)\text{ and }\mathcal{L}%
(^{n}E;F)=\mathcal{L}_{al,2}(^{n}E;F). \label{Tamara}%
\end{equation}

\end{corollary}

Proof. Since every scalar valued $n$-linear ($n\geq2$) mapping defined on
$\mathcal{L}_{\infty}$ spaces is absolutely $(1;2,...,2)$-summing, it is not
hard to prove, using (\ref{DPerez}), that if $E$ is an $\mathcal{L}%
_{\infty,\lambda}$ space, then, regardless of the Banach space $F$, we have
\begin{equation}
\Vert(T(x_{j}^{(1)},...,x_{j}^{(n)})_{j=1}^{m}\Vert_{w,1}\leq\lambda^{n}%
K_{G}3^{\frac{n-2}{2}}\Vert T\Vert\text{ }\Vert(x_{j}^{(1)})_{j=1}^{m}%
\Vert_{w,2}...\Vert(x_{j}^{(n)})_{j=1}^{m}\Vert_{w,2}\label{est}%
\end{equation}
for every continuous $n$-linear mapping $T:E\times...\times E\rightarrow F.$
Then, using the estimates of Proposition \ref{Penelope}, we have
\[
(\int\limits_{0}^{1}\Vert\sum_{j=1}^{m}T(x_{j}^{(1)},...,x_{j}^{(n)}%
)r_{j}(t)\Vert^{2}dt)^{\frac{1}{2}}\leq\Vert(T(x_{j}^{(1)},...,x_{j}%
^{(n)})_{j=1}^{m}\Vert_{w,1}%
\]
and by Definition \ref{esperar} and (\ref{est}), the proof is done.\ The
polynomial case is analogous.$\Box$

\section{A Dvoretzky-Rogers Theorem for almost $p$-summing polynomials}

The Theorem of Dvoretzky-Rogers for absolutely summing linear operators has
natural versions for absolutely summing multilinear mappings and polynomials
(see \cite{Matos1}). A linear Dvoretzky-Rogers Theorem for almost $p$-summing
mappings can be found in \cite[Ex 4.1]{Junek} and tells us that if $p>1$, then
$\mathcal{L}_{al,p}(E;E)\neq\mathcal{L}(E;E)$ for every infinite dimensional
Banach space $E.$ In this section, we will show that we also have multilinear
and polynomial versions for this result.

\begin{lemma}
\label{Beartt}If $P\in\mathcal{P}_{al,p(E)}(^{n}E;F)$ then, regardless of the
$a\in E,$ $dP(a)$ is almost $p$-summing at the origin.
\end{lemma}

Proof. (Adaptation of Lemma 6.1 of\ \cite{Matos1}).We have the following
estimates for $dP(a)(x)$:
\[
dP(a)(x)=\frac{n}{n!2^{n}}\sum\limits_{(e_{i}=1,-1),i=1,...,n}e_{1}%
e_{2}...e_{n}P(e_{1}x+(e_{2}+...+e_{n})a)
\]%
\begin{align*}
&  =\frac{n}{n!2^{n}}\sum\limits_{(e_{i}=1,-1),i=2,...,n}(e_{2}...e_{n}%
P(x+(e_{2}+...+e_{n})a)-(e_{2}...e_{n}P(-x+(e_{2}+...+e_{n})a))\\
&  =\frac{n}{n!2^{n}}(\sum\limits_{(e_{i}=1,-1),i=2,...,n}e_{2}...e_{n}%
[P(x+(e_{2}+...+e_{n})a)-P((e_{2}+...+e_{n})a)])-\\
&  -\frac{n}{n!2^{n}}(\sum\limits_{(e_{i}=1,-1),i=2,...,n}e_{2}...e_{n}%
[P(-x+(e_{2}+...+e_{n})a)-P((e_{2}+...+e_{n})a)])
\end{align*}
Therefore, defining $Q_{e_{2}...e_{n}}(x)=e_{2}...e_{n}[P(x+(e_{2}%
+...+e_{n})a)-P((e_{2}+...+e_{n})a)]$ we have
\[
\int\limits_{0}^{1}\Vert\sum\limits_{j=1}^{k}dP(a)(x_{j})r_{j}(t)\Vert
^{2}dt)^{\frac{1}{2}}=
\]%
\begin{align*}
&  =(\int\limits_{0}^{1}\Vert\sum\limits_{j=1}^{k}\frac{n}{n!2^{n}}%
\sum\limits_{(e_{i}=1,-1),i=2,...,n}(Q_{e_{2}...e_{n}}(x_{j})-Q_{e_{2}%
...e_{n}}(-x_{j}))r_{j}(t)\Vert^{2}dt)^{\frac{1}{2}}\\
&  \leq\frac{n}{n!2^{n}}\sum\limits_{(e_{i}=1,-1),i=2,...,n}(\int
\limits_{0}^{1}\Vert\sum\limits_{j=1}^{k}(Q_{e_{2}...e_{n}}(x_{j}%
)-Q_{e_{2}...e_{n}}(-x_{j}))r_{j}(t)\Vert^{2}dt)^{\frac{1}{2}}\\
&  \leq\frac{n}{n!2^{n}}\{\sum\limits_{(e_{i}=1,-1),i=2,...,n}[(\int
\limits_{0}^{1}\Vert\sum\limits_{j=1}^{k}Q_{e_{2}...e_{n}}(x_{j})r_{j}%
(t)\Vert^{2}dt)^{\frac{1}{2}}+\\
&  +(\int\limits_{0}^{1}\Vert\sum\limits_{j=1}^{k}Q_{e_{2}...e_{n}}%
(-x_{j})r_{j}(t)\Vert^{2}dt)^{\frac{1}{2}}]\}\\
&  \leq\frac{n}{n!2^{n}}\sum\limits_{(e_{i}=1,-1),i=2,...,n}2C_{(e_{2}%
+...+e_{n})a}\Vert(x_{j})_{j=1}^{k}\Vert_{w,p}^{r_{(e_{2}+...+e_{n})a}}\\
&  \leq D\Vert(x_{j})_{j=1}^{k}\Vert_{w,p}^{\min\{r_{(e_{2}...e_{n})a}\}}%
\end{align*}
for $\Vert(x_{j})_{j=1}^{k}\Vert_{w,p}<\delta$ and $0<\delta<\min
\{1;\epsilon_{(e_{2}+...+e_{n})a}\}.$ $\Box$

\begin{theorem}
\label{DR}(Dvoretzky-Rogers for almost $p$-summing polynomials) If $\dim
E<\infty$, then for $p\leq2$ we have
\[
\mathcal{P}_{al,p(E)}(^{n}E;E)=\mathcal{P}(^{n}E;E).
\]
If $\dim E=\infty$ and $p>1$, then $\mathcal{P}_{al,p(E)}(^{n}E;E)\neq
\mathcal{P}(^{n}E;E).$ The multilinear version is also valid.
\end{theorem}

Proof. If $\dim E<\infty$, let us consider $\{e_{1},...,e_{n}\}$ and
$\{\varphi_{1},...,\varphi_{n}\}$ basis for $E$ and $E^{\prime}$ so that
$\varphi_{j}(e_{k})=\delta_{jk}.$ Given an $n$-homogeneous polynomial $P$ from
$E$ into $E$, we have
\[
P(x)=\overset{\vee}{P}(\sum\limits_{j=1}^{m}\varphi_{j}(x)e_{j})^{n}%
=\sum\limits_{j_{1},...,j_{n}=1}^{m}\varphi_{j_{1}}(x)...\varphi_{j_{n}%
}(x)\overset{\vee}{P}(e_{j_{1}},...,e_{j_{n}}).
\]
Since every finite type $n$-homogeneous bounded polynomial is almost
$p$-summing (at zero) for $p\leq2n$ (see \cite[Proposition 3.1 (ii)]{Junek}),
it is not hard to prove that $P$ is almost $p$-summing everywhere$,$ for
$p\leq2$.

On the other hand, suppose that $E$ is an infinite dimensional Banach space.
It suffices to consider the case $1<p\leq2.$ Choose a non null continuous
linear functional $\varphi\in E%
\acute{}%
$ and $a\notin Ker\varphi.$ Define
\[
P(x)=\varphi(x)^{n-1}x.
\]
If we had $P$ almost $p$-summing everywhere, we would have, by Lemma
\ref{Beartt}, $dP(a)$ almost $p$-summing (at zero). Since $\varphi$ is almost
$p$-summing and%
\[
dP(a)(x)=(n-1)\varphi(a)^{n-2}\varphi(x)a+\varphi(a)^{n-1}x,
\]
we would have $\varphi(a)^{n-1}x$ almost $p$-summing. Since $\varphi(a)\neq0$,
we would have that $id_{E}$ is almost $p$-summing, and it is a contradiction.
$\Box$

\begin{example}
It is worth observing that by Corollary \ref{finno}, for $n\geq2$, we have
\[
\mathcal{P}_{al,2}(^{n}c_{0};c_{0})=\mathcal{P}(^{n}c_{0};c_{0})
\]
whereas Theorem \ref{DR} asserts that $\mathcal{P}_{al,2(c_{0})}(^{n}%
c_{0};c_{0})\neq\mathcal{P}(^{n}c_{0};c_{0}).$
\end{example}

Acknowledgment. This paper forms a portion of the author's doctoral thesis,
written under supervision of Professor M\'{a}rio Matos. The author is indebted
to him and to Professor Geraldo Botelho for the suggestions.

\end{document}